\documentclass[10pt]{article}

\usepackage{amsmath}
\usepackage{amsthm}
\usepackage{amssymb}
\usepackage{amsfonts}
\usepackage{graphicx}
\usepackage{verbatim}
\usepackage{mathrsfs}
\usepackage{subfigure}
\usepackage{fancyhdr}
\usepackage{latexsym}
\usepackage{graphicx}
\usepackage{dsfont}

\theoremstyle{plain}
\newtheorem{theorem}{Theorem}[section]

\newtheorem{lemma}[theorem]{Lemma}

\theoremstyle{definition}
\newtheorem{definition}[theorem]{Definition}

\newcommand{\bydef}{\,\stackrel{\mbox{\tiny\textnormal{\raisebox{0ex}[0ex][0ex]{def}}}}{=}\,} 

\newcommand{\R}{\mathbb R}
\newcommand{\C}{\mathbb C}

\newcommand{\hA}{\hat{A}}

\newcommand{\rv}{{\rm v}}
\newcommand{\brv}{\bar{\rm v}}

\newcommand{\bA}{{\boldsymbol{A}}}

\newcommand{\bld}{{\bar{\lambda}}}
\newcommand{\bx}{{\bar x}}
\newcommand{\bl}{{\bar\lambda}}
\newcommand{\cI}{\mathcal{I}}

\title{A method to rigorously enclose eigendecompositions \\ of interval matrices}

\author{Roberto Castelli \thanks{BCAM - Basque Center for Applied Mathematics, Bizkaia Technology Park, 48160 Derio, Bizkaia, SPAIN. 
Email: {\tt rcastelli@bcamath.org}.}
\and
\and Jean-Philippe Lessard\thanks {{\bf Corresponding author}. 
Universit\'e Laval, D\'epartement de Math\'ematiques et de Statistique, Pavillon Alexandre-Vachon, 1045 avenue de la M\'edecine, Qu\'ebec, QC, G1V 0A6, CANADA and BCAM - Basque Center for Applied Mathematics, Bizkaia Technology Park, 48160 Derio, Bizkaia, SPAIN.
Email: {\tt jean-philippe.lessard@mat.ulaval.ca}.
}}

\date{}

\begin{document}

\maketitle

\begin{abstract}
In this paper, a rigorous computational method to enclose eigendecompositions of complex interval matrices is proposed. Each eigenpair $x=(\lambda,\rv)$ is found by solving a nonlinear equation of the form $f(x)=0$ via a contraction argument. The set-up of the method relies on the notion of {\em radii polynomials}, which provide an efficient mean of determining a domain on which the contraction mapping theorem is applicable. 
\end{abstract}

\begin{center}
{\bf \small Keywords} \\ \vspace{.05cm}
{ \small Algebraic eigenvalue problem, Rigorous computations, Contraction mapping theorem}
\end{center}

\begin{center}
{\bf \small Mathematics Subject Classification (2010)} \\ \vspace{.05cm}
{ \small 15A18, 65G40}
\end{center}

\section{Introduction}

Computing eigenvalues and eigenvectors of matrices is a central problem in many fields of applied sciences involving mathematical modelling. When applied to real-life phenomena, any model needs to consider the occurrence of diverse errors in the data, due for instance to inaccuracy of measurements or noise effects. Such uncertainty in the data can be represented by intervals. In the context of studying a matrix with uncertain entries, one can consider an interval matrix $\bA$ whose entries consist of intervals containing all possible errors. Although the entries of $\bA$ are intervals, many classical questions of linear algebra can be raised. For instance, can we demonstrate, given an interval matrix $\bA$, that for any matrix $A \in \bA$, $A$ is invertible, diagonalizable or can we enclose rigorously its eigendecomposition, that is its set of eigenvalues and eigenvectors? In order to address explicitly one of these question, this paper aims at computing rigorous enclosure of eigendecompositions of interval $n\times n$ complex valued matrices $\bA$. More precisely, we propose a method to construct (if possible) a list of balls $\{ B_i : i =1,\dots,n \}$ such that, for every $A \in \bA$ and for each $i = 1,\dots,n$, there exists $x_i=(\lambda_i,\rv_i) \in B_i$ such that $A \rv_i = \lambda_i \rv_i$, where the solution $x_i$ is unique up to a scaling factor of $\rv_i$. 

Before proceeding further, we hasten to mention that several different methods have been developed to find error bounds for computed eigenvalues and eigenvectors of standard (non interval) matrices (e.g. see \cite{MR0566681, MR0684183, MR2068273}). Also, while the problem of computing rigorous bounds for the eigenvalue set of interval matrices is well studied, see for instance \cite{MR1780057, MR2678959} and the references therein, a not so large literature has been produced regarding the simultaneous enclosure of the eigenvalues and eigenvectors of interval matrices. In this direction we refer to \cite{MR0566681, MR1318956}, where different  techniques have been developed to enclose  simple eigenvalues and corresponding eigenvectors, while for double or nearly double eigenvalues a method has been introduced in \cite{MR0843941}. For the rigorous enclosure of multiple or nearly multiple eigenvalues of general complex matrices, a significant contribution has been made by S. Rump in \cite{MR2068273, MR1810531}. 

In this paper, we propose the new idea of enclosing rigorously the eigendecomposition of complex interval matrices by using the notion of {\em radii polynomials}, which provide a computationally efficient way of determining a domain on which the contraction mapping theorem is applicable. The radii polynomials approach, which is very similar to the Krawczyk operator approach, aims at demonstrating existence and local uniqueness of solutions of nonlinear equations. The approach involving the Krawczyk operator consists of applying directly the operator to interval vectors (in the form of small neighbourhoods around a numerical approximation) and then attempt to verify {\em a posteriori} the hypotheses of a contraction mapping argument \cite{MR0255046,MR0657002}. On the other hand, the radii polynomials are {\em a priori} conditions that are derived using analytic estimates, and once they are theoretically constructed, they are used to {\em solve} for the sets (also in the form of small neighbourhoods of a numerical solution) on which a Newton-like operator is a contraction. The advantage of this approach is that most of the estimates are done analytically and generally, and that costly interval arithmetic computations are postponed to the very end of the proofs. It is worth mentioning that the radii polynomials were originally introduced in \cite{MR2338393} to compute equilibria of PDEs with the goal of minimizing the extra computational cost required to prove existence of solutions of infinite dimensional PDEs \cite{MR2776917}. 

In this paper, it is demonstrated that in the context of computing rigorous enclosure of complex matrices, the method based on the radii polynomials is faster than the algorithm introduced in \cite{MR1810531}, with the extra advantage of having local uniqueness (see Section~\ref{sec:comp_cost}). Also, it is demonstrated that an approach based on the Krawczyk operator can be significantly slower than the method involving the radii polynomials (see also Section~\ref{sec:comp_cost}). As in the case of the Krawczyk operator, the radii polynomials verifies existence and uniqueness of a zero of a nonlinear function within the inclusion interval. This also implies multiplicity $1$ of the solution as well as non-singularity of the Jacobian at the solution. Therefore, our new proposed method will necessarily fail to enclose multiple or nearly multiple eigenvalues, a constraint that the method proposed in \cite{MR1810531} does not have.

The paper is organized as follows. In Section~\ref{sec:comp_method}, we introduce the computational method, where we first present the method for non interval matrices. In Section~\ref{sec:interval}, we demonstrate how to generalize the idea to rigorously enclose eigendecompositions of interval matrices. Finally, in Section~\ref{sec:results}, we present applications of our method. In Section~\ref{sec:floquet}, we use the method to rigorously compute the Floquet exponents of a periodic orbit of the Lorenz system of ordinary differential equations. In Section~\ref{sec:large_entries}, we study the applicability of our approach to matrices with interval entries of large radius. Finally, in Section~\ref{sec:comp_cost}, we evaluate the cost of our method and compare it to the cost of the algorithm of S. Rump introduced in \cite{MR1810531} and to a method based on the Krawczyk operator.

\section{The computational method} \label{sec:comp_method}

To begin with, let us fix some notation: throughout this paper we denote by $\mathbb I\C^{n\times n}$ the set of complex matrices with interval entries, $A\in \C^{n\times n}$ an $n\times n$ complex matrix and $\bA\in \mathbb I\C^{n\times n}$ an $n\times n$ interval complex matrix, meaning that any entry of $\bA$ is a complex interval of the form
$$
\bA_{k,j}=[Re(\hat A_{k,j})\pm rad^{(1)}_{k,j}]+i[Im(\hat A_{k,j})\pm rad^{(2)}_{k,j}],  \quad rad^{(1)}_{k,j}, rad^{(2)}_{k,j}\in \R_{+}\ .
$$
The matrix $\hA\in\C^{n\times  n}$ is called the center of $\bA$ while $rad^{(1)}_{k,j}, rad^{(2)}_{k,j}$ are called the radii of the real and imaginary part of $\bA_{k,j}$, respectively.
A matrix $A$ is said to belong to $\bA$, denoted $A\in \bA$, if $A_{k,j}\in \bA_{k,j}$ for any $1 \le k,j \le n$. Bold face letters will always denote interval quantities. Moreover, unless differently specified, 
\begin{itemize}
\item $|\cdot|$ is the complex absolute value and, in case of matrices $M\in\C^{n \times m}$, it acts component-wise, i.e. $|M|_{i,j}=|M_{i,j}|$;
\item given two real matrices $M$, $N$ any relations $<,>,\leq etc$, is assumed component-wise;
\item $I_{n}$ denotes the $n$ dimensional identity matrix, $\mathds 1_{n}$ is the column vector of length $n$ with all the entries equal to 1;
\item given any matrix $M\in \C^{n \times m}$, the object $(M)_{\hat k}$ stands for the $n\times(m-1)$   matrix obtained by deleting the $k$-th column of $M$. 
\end{itemize} 

As already mentioned in the Introduction, given $\bA\in\mathbb I\C^{n\times n}$,  the goal of this paper is to develop a computational method to construct a list of sets $\{B_{i}:i=1,\dots,n\}$ such that, for any $A\in\bA$ and for any $i=1,\dots,n$, there exists  $(\lambda,\rv)\in B_{i}$ solving $A\rv=\lambda\rv$, with $\rv$ unique up to a scaling factor.

To simplify the exposition, we first present the method in the context of non interval matrices $A\in\C^{n\times n}$, that is one introduces a method to enclose the solutions $(\lambda,\rv)$ of the equation
\begin{equation}\label{eq:eig}
A\rv=\lambda\rv .
\end{equation}
As one shall see in Section~\ref{sec:interval}, only minor modifications are necessary for the extension to the interval case.
 
Suppose that an approximate eigenpair of $A$ has been computed, that is $(\bld,\brv)$ such that $A \brv \approx \bld \brv$ and let $f(x)$ be the function $f:\C^{n}\rightarrow \C^n$ that maps a point $x=(\lambda,\rv_1,\rv_2,\dots,\rv_{k-1},\rv_{k+1},\dots,\rv_n)$ to 
\begin{equation}\label{eq:f}
f(x)=A\left[\begin{array}{c}
\rv_{1}\\
\vdots\\
\brv_{k}\\
\vdots\\
\rv_{n}
\end{array}\right]-\lambda\left[\begin{array}{c}
\rv_{1}\\
\vdots\\
\brv_{k}\\
\vdots\\
\rv_{n}
\end{array}\right]
\end{equation}
where $\brv_{k}$ is the largest component of $\brv$.  By definition, a solution $x$ of  $f(x)=0$ corresponds to an eigenpair $(\lambda,\rv)$  of $A$ with the eigenvalue $\lambda$ given by the first component of $x$ and the eigenvector  $\rv=(\rv_{1},\dots,\rv_{k-1},\brv_{k},\rv_{k+1},\dots,\rv_{n})$. 
Thus, the target of the incoming analysis is to prove the existence and to provide rigorous bounds for the zeros of $f(x)$. Note that the unknowns in the equation $f(x)=0$ are $\lambda$ and $n-1$ component of $\rv$, while the remaining component $\brv_{k}$ is a fixed parameter of the problem. Since the  eigenvectors are invariant under rescaling, the solutions of \eqref{eq:eig} come in continuous families. However, fixing one of the component of $\rv$ (in our case letting $\rv_k=\brv_k$) removes such arbitrariness and therefore isolates the zeros of $f$. Denoting $\bar x=(\bld,\brv_{1},\brv_2,\dots,\brv_{k-1},\brv_{k+1},\dots,\brv_n) $ and $Df(\bar x)$ the Jacobian matrix of $f$ at $\bar x$, one has that
\begin{equation}\label{eq:jac}
Df(\bar x)=\left[ 
\left[\begin{array}{c}
\brv_{1}\\
\vdots\\
\brv_{k}\\
\vdots\\
\brv_{n}
\end{array}\right](A-\bld I_{n})_{\hat k}
\right].
\end{equation}
The problem of finding the zeros of the function $f(x)$ is addressed by introducing an operator $T$ on a Banach space whose fixed points correspond to solutions of  $f(x)=0$. 

Let $\Omega$ be the Banach space 
$$
\Omega=\{x\in\C^{n}: ||x||_{\Omega}<\infty\},\quad ||x||_{\Omega}=\max_{i}\{|x_{i}| \},
$$
and define the operator $T:\Omega\rightarrow \Omega$ by
\begin{equation}\label{eq:T}
T(x)=x-Rf(x),
\end{equation}
where $R$ is a numerical  inverse of $Df(\bar x)$, i.e. $R \cdot Df(\bar x)\approx I_{n}$.
Since fixed points of $T$ correspond to zeros of $f(x)$, the idea is to construct a small set $B \subset \Omega$ such that $T:B \rightarrow B$ is a contraction, and then to apply the contraction mapping theorem to conclude about the existence of a unique fixed point of $T$ in $B$. Note that $\bar x$ is an approximate zero of $f$ and  the operator $T$ has been defined as a Newton-like operator around the point $\bar x$, thus it is advantageous to test the contractibility of $T$ on neighbourhoods of $\bar x$ in $\Omega$.

More precisely, denote by
$$
B(r)=\{x\in\Omega, ||x||_{\Omega}\leq r\}
$$
the close ball of radius $r$ around the origin and let $B_{\bar x}(r)=\bar x+ B(r)$ the  ball with the same radius and  centered at $\bar x$. Treating $r$ as a variable, we choose the balls $B_{\bar x}(r)$ as the candidate sets where to check if $T$ is a contraction.
The question whether $T$ is a contracting map will be formulated in terms of the verification of a set of computable conditions, called {\it radii polynomials}, whose construction is based on verifying efficiently the hypothesis of the next result.

\begin{theorem}\label{th:radii}
Suppose that $Y,Z(r)\in \Omega$ are such that  
\begin{equation}\label{eq:hyYZ}
|T(\bx)-\bx|\leq Y,\quad \sup_{b,c\in B(r)}|DT(\bx+b)c|\leq Z(r),
\end{equation}
and satisfy
$$
\|Y+Z(r)\|_{\Omega}<r.
$$
Then there exists a unique $x\in B_{\bar x}(r)$ such that $T(x)=x$.
\end{theorem} 
\proof
The mean value theorem applied component-wise to $T$ implies that for any $x,y\in B_{\bar x}(r)$ and for any $k=1,\dots,n$,
$$ 
T_{k}(x)-T_{k}(y)=DT_{k}(z)(x-y),\quad {\rm for}\ z\in\{tx+(1-t)y: t\in[0,1] \}\subset B_{\bar x}(r).
$$
Then, 
\begin{equation}\label{eq:bnd}
|T_{k}(x)-T_{k}(y)|=\left|DT_{k}(z)\frac{r(x-y)}{\|x-y\|_{\Omega}}\right|\frac{1}{r}\|x-y\|_{\Omega}\leq \frac{Z_{k}(r)}{r}\|x-y\|_{\Omega}.
\end{equation}
Choosing $y=\bar x$ and using the triangular inequality, one has that
$$
|T_{k}(x)-\bar x_{k}|\leq |T_{k}(x)-T_{k}(\bx)|+|T_{k}(\bx)-\bx_{k}|\leq Y_{k}+Z_{k}(r)\leq r,
$$
where the last inequality follows from the fact that $Y,Z\in \R^{n}_{+}$, and thus $\|Y+Z(r)\|_{\Omega}=\max_{i}\{Y_{i}+Z_{i}(r)\}$. That proves that  $T(B_{\bx}(r))\subseteq B_{\bx}(r)$. From \eqref{eq:bnd}, it follows that
$$
\|T(x)-T(y)\|_{\Omega}=\max_{k}\{|T_{k}(x)-T_{k}(y)|\}\leq\frac{\|Z(r)\|_{\Omega}}{r}\|x-y\|_{\Omega}.
$$
Since $\|Z(r)\|_{\Omega}<r$, $T$ is a contraction on $B_{\bar x}(r)$. Thus, from the contraction mapping theorem, the exists a unique fixed point of $T$ in $B_{\bar x}(r)$.

\qed

\begin{definition} \label{def:rad_poly}
Given the vectors $Y,Z(r)\in \Omega$ satisfying \eqref{eq:hyYZ}, we define the {\em radii polynomials} $p_k(r)$, $k=1,\dots, n$ by
$$
 p_k(r)=(Y+Z(r))_k-r .
$$
\end{definition}
The following result holds.
\begin{lemma}\label{lemma:uniq}
Consider the radii polynomials $p_k(r)$, $k=1,\dots, n$ introduced in Definition~\ref{def:rad_poly}. Then for any $r>0$ such that 
$$
p_k(r)<0, \quad {\rm ~Êfor~Êall~} k=1,\dots,n,
$$
there exists a unique $x\in B_{\bx}(r)$ such that  $f(x)=0$.
\end{lemma}
\proof
Suppose $\bar r>0$ is such that $p_k(\bar r)<0, \quad {\rm ~Êfor~Êall~} k=1,\dots,n$. 
  For  \eqref{eq:hyYZ} all the entries of  $Y,Z(\bar r)$ are real positive,  thus   $\max_{i}\{(Y+Z(\bar r))_{i}\}=||Y+Z(\bar r)||_{\Omega}<\bar r$. From Theorem~\ref{th:radii} there exists a unique $x\in B_{\bx}(\bar r)$ such that $T(x)=x$ and therefore $f(x)=0$.
  
   \qed

\noindent We proceed with the explicit construction of the bounds $Y$ and $Z$.  Since $T(\bx)-\bx=-Rf(\bx)$, let 
\begin{equation}\label{eq:Y}
Y \bydef |Rf(\bx)|.
\end{equation}
The bound $Z(r)$ satisfying \eqref{eq:hyYZ} is constructed as a polynomial in $r$. First rewrite $DT(\bx+b)c$ as
$$
DT(\bx+b)c= (I_{n}-R \cdot Df(\bx))c+R[(Df(\bx)-Df(\bx+b))c]\\
$$
so that
\begin{equation}\label{eq:split}
|DT(\bx+b)c|\leq |(I_{n}-R \cdot Df(\bx))c|+|R[(Df(\bx)-Df(\bx+b))c]|.
\end{equation}
Define
\begin{equation}\label{eq:Z}
Z(r)=rZ_{0}+r^{2}Z_{1},
\end{equation}
where
\begin{equation}\label{eq:Z0}
Z_{0}\bydef|I_{n}-R \cdot Df(\bx)|\mathds 1_{n},\quad
Z_{1}\bydef 2|R| \hat{\mathds 1}_{n} ,
\end{equation}
where $\hat{\mathds 1}_{n}$ is the same as $\mathds 1_{n}$ with a zero, instead of one,  in the $k$-th component.   
\begin{lemma}
Consider the polynomial vector $Z(r)$ defined by \eqref{eq:Z}. Then
$$
Z(r)\geq  \sup_{b,c\in B(r)}|DT(\bx+b)c| .
$$
\end{lemma}
\proof
From \eqref{eq:split}, the statement follows by proving  that
\begin{itemize}
\item[i)]  ${\displaystyle \sup_{c\in B(r)}|(I_{n}-R \cdot Df(\bx))c|\leq rZ_{0}}$
\item[ii)] ${\displaystyle  \sup_{b,c\in B(r)}|R[(Df(\bx)-Df(\bx+b))c]|\leq r^{2}Z_{1} }$.
\end{itemize}
 Since $c\in B(r)$, $|c|\leq r\mathds 1_{n}$, it follows that
 $$
 \sup_{c\in B(r)}|(I_{n}-R \cdot Df(\bx))c|\leq rZ_{0}\ .
$$
That proves $i)$. For any $b=(b_{\lambda},b_{1},\dots, b_{k-1},b_{k},\dots, b_{n})$
\[
(Df(\bx)-Df(\bx+b))=
\left[-\left[\begin{array}{c}
b_{1}\\
\vdots\\
b_{k-1}\\
0\\
b_{k+1}\\
\vdots\\
b_{n}
\end{array}\right](b_{\lambda} I_n)_{\hat k}
\right].
\]
Note that the $k$-th row of the above matrix is null. Since $|b_{i}|\leq r$, we have $|(Df(\bx)-Df(\bx+b))c|\leq 2r^{2}\hat{\mathds 1}_{n}$ and therefore
 $\sup_{b,c\in B(r)}|R[(Df(\bx)-Df(\bx+b))c]|\leq 2r^{2}|R|\hat{\mathds 1}_{n}=r^{2}Z_{1} $.
 
\qed

In summary, given an approximate eigenpair $(\bar\lambda,\brv)$,  the method consists of computing rigorously the bounds $Y$ and $Z(r)$ given by \eqref{eq:Y} and \eqref{eq:Z}, and then to check whether there exists an interval $\cI$ where all the polynomials  $p_k(r)$ are negative. If $\cI \neq \emptyset$ we select   $r=\inf \cI$ and we conclude that $f(x)=0$ has a unique  solution within the ball $B_{\bar x}(r)$. In practice, we get the existence of an eigenpair $(\lambda,\rv)$ of $A$, with $|\lambda-\bar\lambda|\leq r$, $|\rv_{j}-\bar \rv_{j}|\leq r$, for $j\neq k$ and $\rv_{k}=\bar \rv_{k}$. To prove the existence of a second eigenpair of $A$, it is necessary to provide an approximate solution  $(\bar\lambda,\brv)$, different from the previous one, and to repeat the computation.  

\subsection{Extension to the interval case}\label{sec:interval}

Besides few modifications necessary to deal with interval quantities,  the procedure to compute rigorously bounds for the eigendecomposition of an interval matrix $\bA\in \mathbb I\C^{n\times n}$ is basically the same as for the scalar case. However, a fundamental difference consists in the fact that all the computations are done in   the interval arithmetic regime \cite{MR1780057}, in which any of  the basic operations $\circ\in\{+,-,\cdot,/\}$  is extended to the interval case in order to satisfy the general assumption
\begin{equation}\label{interval}
\forall P\in\boldsymbol P\quad \forall Q\in\boldsymbol Q, \quad P\circ Q\in \boldsymbol P \circ \boldsymbol Q\ .
\end{equation}

Given $\bA$  an interval complex valued matrix, we now address the problem we stated at the beginning of the paper, that is how to construct  sets $\{B_{i}\}_{i=1}^{n}$ so that each $A\in\bA$ admits one and only one eigenpair $(\lambda,\rv)$ in any of the $B_{i}$'s.  
Recall that $\hA$ is the center of the interval matrix $\bA$. We first compute $(\bl, \brv)$ an approximate eigenpair of $\hA$ and, as before, define $\bar x=(\bld,\brv_{1},\brv_2,\dots,\brv_{k-1},\brv_{k+1},\dots,\brv_n) $ where the missing component is chosen so that $\brv_{k}=\max_{j}\{\brv_{j}\}$. Then, replacing the scalar matrix $A$ in \eqref{eq:f} by the interval matrix $\bA$, the function $f(x)$ and the Jacobian matrix $Df(\bar x)$ defined in \eqref{eq:f} and \eqref{eq:jac} are replaced respectively by $\boldsymbol{f}:\C^{n}\rightarrow \mathbb I\C^{n}$  and by an interval matrix $\boldsymbol{Df}(\bar x)$ that represents a linear operator from $\C^{n}$ to $\mathbb I\C^{n}$. We choose $ R$ to be a numerical inverse of $\widehat{Df}(\bar x)$, the center of $\boldsymbol{Df}(\bar x)$, and we proceed to the definition of the operator $\boldsymbol{T}(x)=x-R\boldsymbol{f}(x)$ and to  the bounds $Y$ and $Z(r)$, as done before with the boldface quantities in place of the previous one.   Clearly the quantities on the left hand side of relations \eqref{eq:hyYZ} are now intervals, thus we  define component-wise $Y$, $Z_{1}$, $Z_{2}$  as the supremum of the intervals appearing on the right hand sides of \eqref{eq:Y} and \eqref{eq:Z0}. That will yield the uniform bounds
$$
| \boldsymbol T(\bx)-\bx|\leq Y,\quad \sup_{b,c\in B(r)}|D\boldsymbol T(\bx+b)c|\leq Z(r),
$$
where $D\boldsymbol T(x) = I_n - R \cdot D \boldsymbol f(x)$.

Suppose that $r>0$ is such that the radii polynomials $p_k(r)<0$ for all $k$. Then interval arithmetic insures that forÊ all $A\in\bA$, there exists a unique $(\lambda,\rv)$ such that $|\lambda-\bar\lambda|\leq r$, $|\rv_{j}-\bar \rv_{j}|\leq r$, $\rv_{k}=\bar \rv_{k}$, and $A\rv=\lambda \rv$. In other words, $r$ is a uniform bound in $\bA$ for the existence of an eigenpair of any $A\in\bA$. Indeed, having fixed $(\bl,\brv)$ and $ R\approx \big(\widehat{Df}(\bar x)\big)^{-1}$, for any $A\in\bA$ define $f_{A}(x)$ and $Df_{A}(\bx)$ as in \eqref{eq:f} and \eqref{eq:jac}, and the fixed point operator $T_{A}(x)=x-Rf_{A}(x)$. The fundamental inclusion \eqref{interval} implies  that   $f_{A}(x)\in \boldsymbol{f}(x)$, $Df_{A}(\bar x)\in \boldsymbol{Df}(\bar x)$ and $T_{A}(x)\in\boldsymbol{T}(x)$, for any $A\in\bA$ and $x\in\C^{n}$. Thus, as $A$ varies in $\bA$, the bounds \eqref{eq:hyYZ}, with $T_{A}$ in place of $T$,  are satisfied for the same $Y, Z, r$ proving the existence of a fixed point in $B_{\bar x}(r)$ for any $T_{A}$ and consequently an eigenpair for any $A\in \bA$.

\section{Results} \label{sec:results}

In this section we report some computational results. All the computations have been done in {\it Matlab} supported by the package {\it Intlab} \cite{Ru99a} where the interval arithmetic routines have been implemented.  The approximate eigenpairs $(\bar \lambda,\bar\rv)$ of $\hA$ 
 have been computed running the standard \verb#eig.m# function in {\it Matlab}.  In order to avoid rounding error and to obtain rigorous results, we emphasize that   the computational algorithm treats any matrix as an interval matrix. Thus, even if one wishes to deal with a scalar matrix $A$, the method first constructs  a (narrow) interval matrix around $A$ and perform all the computation with interval arithmetics. 

 \subsection{Example 1: rigorous computations of Floquet exponents} \label{sec:floquet}
The first example concerns the rigorous enclosure of the Floquet exponents and related eigenvectors associated to a periodic orbit of the Lorenz system
\begin{equation}\label{syst:lor}
\left\{
\begin{array}{l}
\dot u_{1}=\sigma (u_{2}-u_{1})\\
\dot u_{2}=\rho u_{1}-u_{2}-u_{1}u_{3}\\ 
\dot u_{3}=u_{1}u_{2}-\beta u_{3}
\end{array} \right. 
.
\end{equation}
For a choice of the parameters $\sigma=10$, $ \beta=8/3$, $\rho=20.8815$ the existence of a periodic orbit $\gamma(t)$ in a neighborhood of an approximate solution has been proved \cite{CGL}. It is known that the stability character of the orbit $\gamma(t)$ is encoded by the Floquet exponents, which are the eigenvalues of a particular real matrix, denoted by $A$,  resulting by integrating the linearized system around $\gamma(t)$. Without going into the details (we refer to \cite{floquet} for an exhaustive explanation), we only mention that one of the Floquet exponents is zero, due to the  time shift, while the number of the other eigenvalues of $A$ with negative (positive)  real part gives the dimension of the stable (unstable) manifold.   Moreover the associated eigenvectors provide the directions tangent to the invariant manifolds on $\gamma(t)$ (tangent bundles).

By means of a computational method based on the radii polynomials, in \cite{floquet} we proved that the matrix $A$ associated to the solution $\gamma(t)$ lies within the interval matrix  $\bA$ centered at 
$$\hat A=\left[
\begin{array}{rrr}
-10.55360193   &5.33379647 &-5.24740415\\
   0.31403414  &2.33062549  &-3.32865541\\
  -7.49045333   &5.01386821  &-5.44369022\\
\end{array}\right]
$$
with radius $rad= 9.66146973\cdot10^{-7}$, meaning that each entry $\bA(i,j)$ consists of the interval $[\hat A(i,j)-rad,\hat A(i,j)+rad]$.
 Following the computational method discussed in Section~\ref{sec:comp_method}, we compute the enclosure of the eigenpairs of $\bA$: it results that any $A \in \bA$ admits three eigenpairs $(\lambda_{i},\rv_{i})$, $i=1,2,3$ each one lying in the ball of radius $r_{i}$ around the approximate values $(\bar{\lambda}_{i},\bar{\rv}_{i})$ given in Table~\ref{tab:floquet_data}.
 
 \begin{table}[htdp]
 $$
 \begin{array}{c|c|c|c}
 & i=1&i=2&i=3\\
 \hline
 \vspace{-5pt}&&&\\
 r_{i}&    2.774764083439355\cdot 10^{-6}&  3.567796353801448\cdot 10^{-5} & 3.649406638638561\cdot 10^{-5}\\ 
 \hline
 &&& \\
 \bar{\lambda}_{i}& -13.962049368058944 & -9.363255645359691\cdot 10^{-14}& 0.295382701392333\\
 &&&\\
 \bar{\rv}_{i}&\left[\begin{array}{l}
   -0.777880129851759\\
  -0.111361799590875\\
  -0.618466695282529\\
 \end{array}\right]&
 \left[\begin{array}{l}
 0.121332037790074\\
   0.807699961688802\\
   0.576974270218020\\
    \end{array}\right]&
    \left[\begin{array}{l}
    0.156626754180928\\
   0.835985628946308\\
   0.525924032603563\\
      \end{array}\right]
 \end{array}
 $$
 \vspace{-.5cm}
 \caption{Data associated to the rigorous computation of the Floquet exponents of the periodic orbit $\gamma$.}
 \label{tab:floquet_data}
 \end{table}

We remark that in the general situation the genuine solution $(\lambda,\rv)$ of the eigenproblem is proved to exist in a complex neighborhood of the approximate solution $(\bar\lambda,\bar\rv)$.  Therefore, even if one or both $\bar\lambda$ and $\bar\rv$ are real vectors, the same can not be concluded for $\lambda$ or $\rv$. However, if the  matrix $A$ and the approximate solution $\bar\lambda$ and $\bar\rv$ are real and the computation is successful, then the genuine solution so obtained by solving the radii polynomials is also real. Indeed, suppose  the contrary, that is the exact solution $\lambda$ and $\rv$ are complex. Since $A$ is real,  the complex conjugate couple $(\mathcal C(\lambda), \mathcal C(\rv)))$ is also a solution of the eigenproblem, $A\mathcal C(\rv)=\mathcal C(\lambda)\mathcal C(\rv)$. But both the solutions $(  \lambda,\rv)$ and $(\mathcal C(\lambda), \mathcal C(\rv)))$ belong to the same ball in $\Omega$ around $\bar x$ and this violates the uniqueness result stated in Lemma \ref{lemma:uniq}.  The same argument extends in the case of interval matrices.

Coming back to the previous example, we conclude that  the  Floquet multipliers of $\gamma(t)$ are real, one is negative and one is positive. Note that the zero Floquet multiplier is indeed contained in the ball around $(\bar\lambda_{2},\bar{\rv}_{2})$.
 
 \subsection{Example 2: matrices with interval entries of large radius} \label{sec:large_entries}
 
\begin{figure}
\begin{center}
 \includegraphics[scale=0.6]{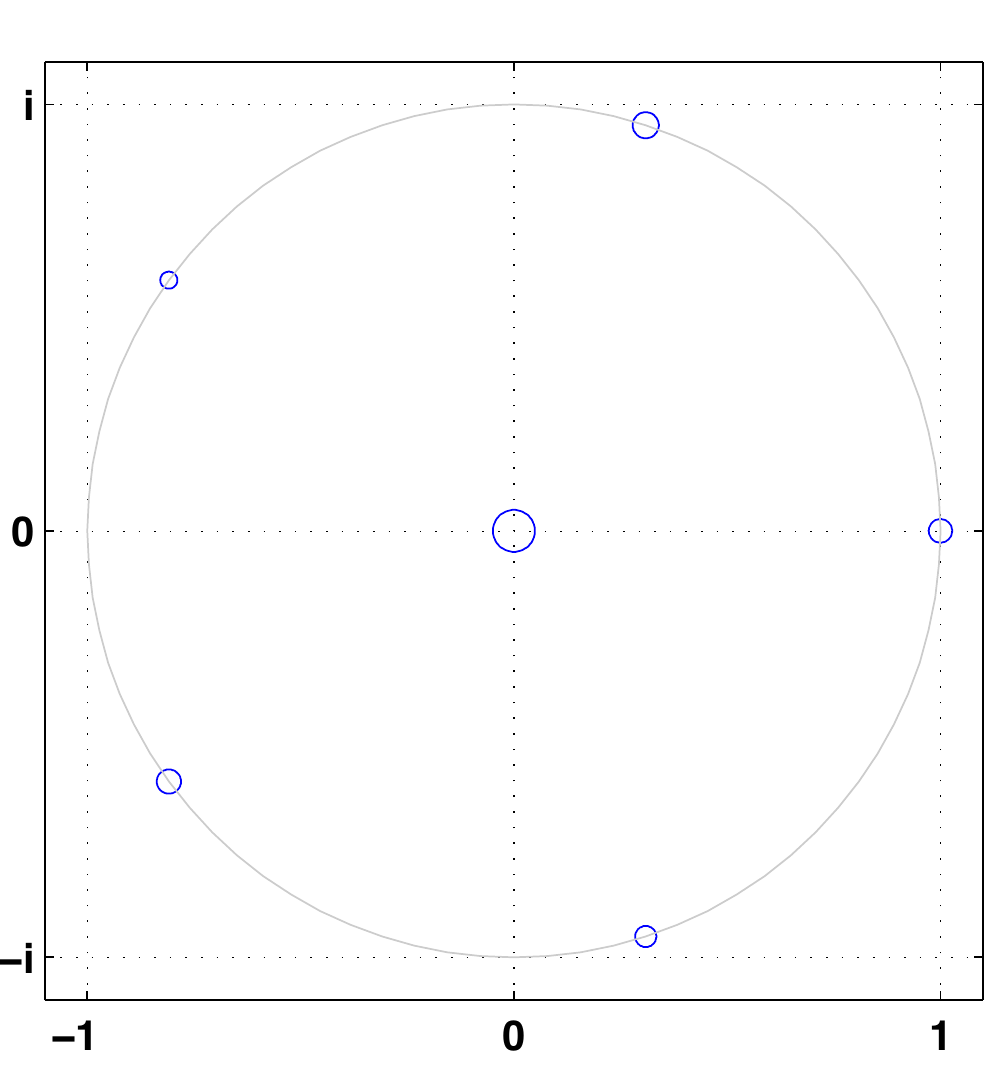}
 \end{center}
 \vspace{-.3cm}
 \caption {\small Balls in the complex plane enclosing the 6 eigenvalues of any $A\in\bA$ for $rad=1.3\cdot 10^{-3}$.}
\label{figure}
\end{figure}

In the next sample we compute the eigendecomposition of an interval matrix $\bA$ constructed as follows: consider the complex number $\lambda_{0}=0$ and $\lambda_{j}=e^{i\frac{2\pi}{5}j}$, $j=1,\dots,5$ and define $D$ as the diagonal matrix with entries $\lambda_{k}$, $k=0,\dots,5$. Let $\hA=XDX^{-1}$, for a random matrix $X$ with values in the complex square $[-1,1]+i[-1,1]$ and finally let $\bA$ be the interval complex matrix centered at $\hA$ with component-wise radius $rad$ both in the real and imaginary part. For different values of $rad$ we compute the  enclosure of the eigenvalues of $\bA$: given  $\bar\lambda_{k}$ the approximate eigenvalues of $\hA$,
$$
\begin{array}{rr}
\bl_{0}=&0.00000 + 0.00000i\\
\bl_{1}=& 0.30901 + 0.95105i\\
\bl_{2}=&-0.80901+ 0.58778i\\
\bl_{3}=&-0.80901 - 0.58778i \\
\bl_{4}=&0.30901 - 0.95105i\\
\bl_{5}=&1.00000-0.00000i
\end{array}
$$
denote by $r_{k }$, $k=0,\dots 5$  the radius of the ball in the complex plane centered at $\bar\lambda_{k}$ inside which, for any $A\in\bA$,  a unique eigenvalues of $A$ has been proved to exist. The results are presented in Table \ref{Tab:1}, where the different values of $rad$ are given in the first column,  while each row collects the result for the radii $r_{k}$.   For values of $rad\leq 1.3\cdot 10^{-3}$ the method succeeded in computing the enclosure of  the entire eigendecomposition of $\bA$, while for larger  values of $rad$ the method starts failing in enclosing some of the eigenpairs, up to the value $rad=3.2\cdot 10^{-3}$, where no computation is successful. Figure \ref{figure} shows the radii of the disks in the complex plane enclosing the six eigenvalues for the case $rad=1.3\cdot 10^{-3}$.

\begin{table}[htdp]
\begin{center}
\begin{tabular}{c|c|c|c|c|c|c}
$rad$&$r_{0}$&$r_{1}$&$r_{2}$&$r_{3}$&$r_{4}$&$r_{5}$\\
\hline&&&&&\\
$1\cdot 10^{-5}$&$2.0\cdot 10^{-4}$&$1.5\cdot 10^{-4}$&$1.6\cdot 10^{-4}$&$1.3\cdot 10^{-4}$&$1.8\cdot 10^{-4}$&$1.6\cdot 10^{-4}$\\
$1\cdot 10^{-4}$&$0.0021$&$0.0016$&$0.0018$&$0.0014$&$0.0019$&$0.0017$\\
$1\cdot 10^{-3}$&$0.0281$&$0.0181$&$0.0202$&$0.0149$&$0.0219$&$0.0197$\\
$1.3\cdot 10^{-3}$&$0.0492$&$0.0247$&$0.0284$&$0.0201$&$0.0306$&$0.0276$\\
$1.9\cdot 10^{-3}$& $-$ &$0.0416$&$0.0563$&$0.0323$&$0.0515$&$0.0581$\\
$2\cdot 10^{-3}$& $-$ &$0.0452$& $-$ &$0.0346$&$0.0678$&$0.0591$\\
$2.5\cdot 10^{-3}$& $-$ &$0.0759$& $-$ &$0.0485$& $-$ & $-$ \\
$3.1\cdot 10^{-3}$& $-$ & $-$ & $-$ &$0.0828$& $-$ & $-$ \\
$3.2\cdot 10^{-3}$& $-$ & $-$ & $-$ & $-$ & $-$ & $-$ \\
\end{tabular}
\caption{Enclosures of the eigenvalues of a complex interval matrix $\bA$, as the radius $rad$ of the entries of $\bA$ increases. }
\label{Tab:1}
\end{center}
\end{table}

\subsection{Example 3: evaluating the performance} \label{sec:comp_cost}

In this last section we compare the performance of our new method, that we denote here by $radiipol$, with two different algorithms developed by S. Rump. The first one, here denoted by {\it verifyeig}, has been introduced in  \cite{MR1810531} with the primary goal of computing enclosures of multiple of nearly multiple eigenvalues (and related eigenvectors) of interval matrices. It consists of a verification method and provides rigorous bounds around an approximate solution within which the existence, but not uniqueness, of the exact solution of the eigenproblem is proved. The second one, here denoted by {\it verifynlss}, is based on a Krawczyk operator \cite{MR0255046,MR0657002} and is a general routine to rigorously compute well separated zeros of nonlinear functions. In fact, in the code \verb#verifyeig.m# (available in the library {\em Intlab} \cite{Ru99a}), where the method {\it verifyeig} has been implemented, the author suggested to use  {\it verifynlss} to compute simple and well separated eigenpairs.  This method is implemented in the code \verb#verifynlss.m# in the library {\em Intlab} \cite{Ru99a}.

We summarize the obtained results in Table~\ref{Tab:time} and Table~\ref{Tab:performance}. In Table~\ref{Tab:time}, the computational time necessary to compute the entire eigendecomposition of a test matrix (measured in seconds on a 2.4 GHz computer using the \verb#tic-toc# Matlab function) is reported, independently whether the methods were successful or not. Table~\ref{Tab:performance} presents the output of the algorithms: for each computation of the eigendecomposition, it provides the average of the radius of the balls enclosing the exact eigenpairs. The entry $-$ means that the method fails in the enclosure of at least one of the eigenpair.

For both experiments the test matrices $\bA$ have been constructed as in the previous section: given $N$ we define $D\in \C^{N+1,N+1}$ as a diagonal matrix with entries given by $N$ equispaced values on the unit circle in the complex plane and $0$, i.e. $diag(D)=[0,e^{i\frac{2\pi}{N} j}], j=1,\dots, N$. Then let  $\hA=XDX^{-1}$, where $X$ is a complex random matrix with entries in the complex square $[-1,1]+i[-1,1]$ and finally define $\bA$ as the interval complex matrix centered in $\hA$ and of radius $rad$.

From the results of Table \ref{Tab:time}, we conclude that our new proposed method $radiipol$ is much faster than the method $verifynlss$ based on the Krawczyk operator, while also proving existence and uniqueness of the solutions. Moreover, our algorithm is faster than $verifyeig$, especially for small $N$, and the computational time of the two methods is almost equivalent for large $N$. This achievement have been possible thanks to the analytical estimates (in the form of the radii polynomials) introduced in the method, that allow minimizing the number of computations done with interval quantities. Hence, the $radiipol$ method provides a computationally efficient way to compute eigendecompositions of complex interval matrices.

The results presented in Table~\ref{Tab:performance} confirm that the new approach $radiipol$ is satisfactory also from the point of view of the accuracy of the results. Indeed, while the algorithm $verifynlss$ fails quite soon as $N$ and $rad$ increase (it fails for $rad=0$ and for all $N\geq 15$), the new algorithm is successful also for large entries of $\bA$. Moreover, as one can see in Table~\ref{Tab:performance}, the performance of $radiipol$ is very close to the performance of the algorithm $verifyeig$ which, we underline this one last time, does not prove uniqueness of the solution.

\begin{table}
\begin{minipage}[c]{0.5\linewidth}\centering
\qquad Computational time (s)\\
\begin{tabular}{c|c|c|c}\hline&&\\
$N$&$radiipol$& $verifyeig$ & $verifynlss$\\
\hline&&\\
5&0.0148&0.0372&0.3423\\
10&0.0282&0.0650&0.6432\\
50&0.2481&0.4170& 7.4476${}^*$ \\
100&1.4580&1.8814&66.119${}^*$\\
200&14.512&16.739& 792.76${}^*$ \\
500& 377.3 &423.6 &
\end{tabular}
\end{minipage}
\begin{minipage}[c]{0.5\linewidth}\centering
\vspace{10pt}
\includegraphics[scale=0.35]{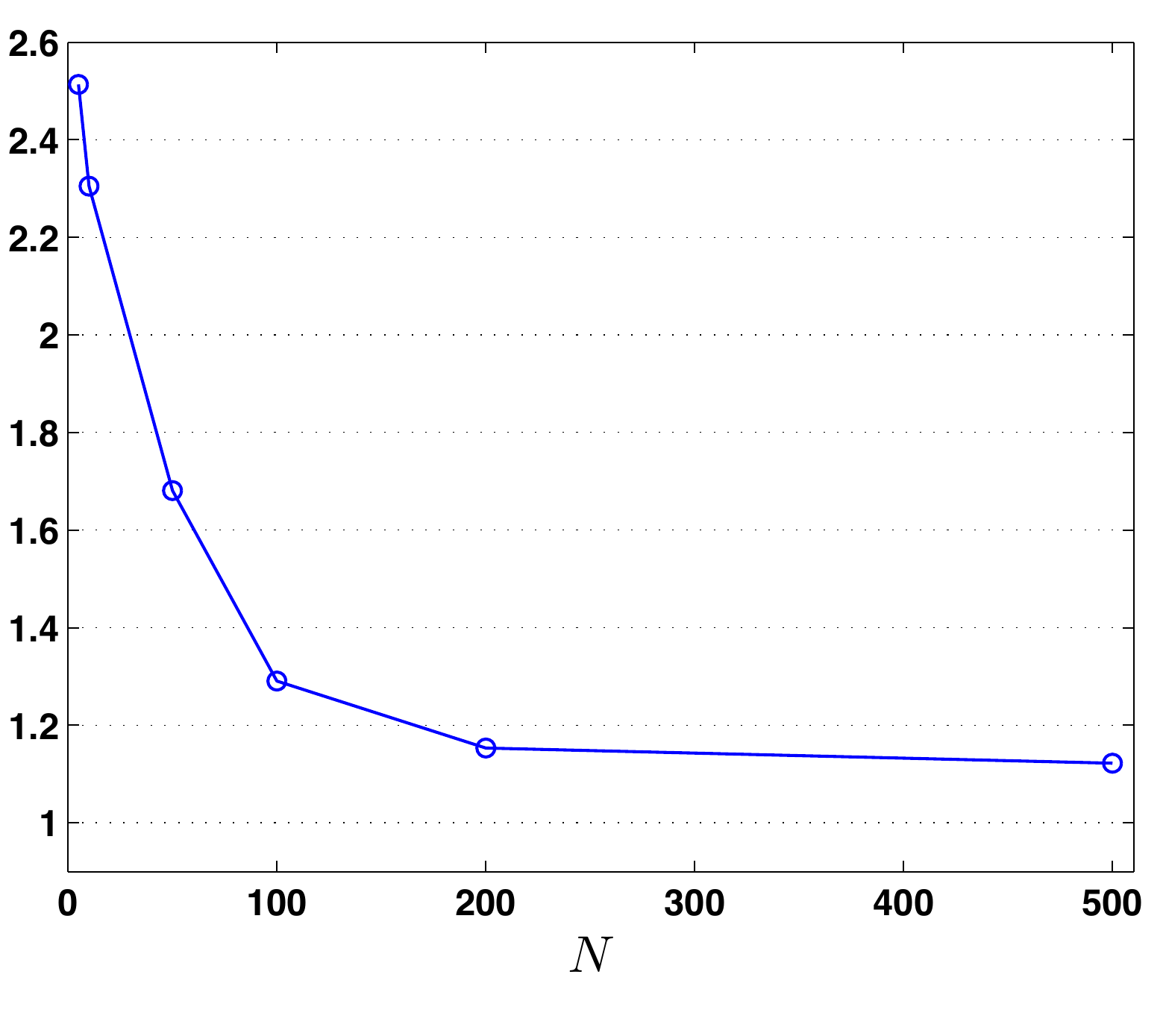}
\end{minipage}
\vspace{-.6cm}
\caption{\small Comparison of the running time necessary to enclose the entire eigendecomposition of $\bA$ using $radiipol$ and the two other methods $verifyeig$ and $verifynlss$. $\bA$ is a complex interval matrix of dimension $N+1$ and $rad=10^{-15}$. The entries of the form ${(\cdot)}^*$ correspond to unsuccessful computations. The figure on the right depicts the ratio between the computational time of the algorithms $verifyeig$ and $radiipol$.} 
\label{Tab:time}
\end{table}

\begin{table}
\begin{center}
\qquad \qquad Average of the radius of the disks enclosing the eigenvalues\\
\end{center}
\begin{tabular}{c|ccc|ccccc}
&& N=5&&&N=10&&\\
\hline
$rad$&$10^{-20}$ &$10^{-10}$&$10^{-4}$&$10^{-10}$&$10^{-5}$&$10^{-4}$&$10^{-3}$\\
\hline 
$radiipol$& 9.14$\cdot 10^{-15}$& 2.76$\cdot 10^{-9}$&0.0019&3.25$\cdot 10^{-9}$&4.61$\cdot 10^{-4}$& $-$&$-$\\
$verifyeig$&4.69 $\cdot 10^{-15}$& 2.07$\cdot 10^{-9}$&0.0016&2.08$\cdot 10^{-9}$&3.02$\cdot 10^{-4}$&0.0049&$-$\\
$verifynlss$&6.26 $\cdot 10^{-9}$& $-$& $-$&$-$&$-$&$-$&$-$\\
\end{tabular}

\vspace{5pt}
\begin{tabular}{c|ccc|ccc|c}
&& N=50&&&N=100&& N=150\\
\hline
$rad$&$10^{-10}$ &$10^{-8}$&$10^{-5}$&$10^{-10}$&$10^{-8}$&$10^{-7}$& $10^{-10}$\\
\hline 
$radiipol$ & 2.69$\cdot 10^{-7}$& 4.94$\cdot 10^{-5}$&$-$&9.02$\cdot 10^{-7}$&$-$&$-$& 1.31$\cdot 10^{-6}$\\
$verifyeig$\ &5.59$\cdot 10^{-8}$& 9.45$\cdot 10^{-6}$&$-$&1.31$\cdot 10^{-7}$&2.07$\cdot 10^{-5}$&$-$&1.64$\cdot 10^{-7}$\\
$verifynlss$& $-$& $-$&$-$&$-$&$-$&$-$& $-$
\end{tabular}
\caption{\small Comparison of the accuracy  of the three methods  as the dimension $N$ and the radius $rad$ of the test matrix $\bA$ change. The entry $-$ means that the method fails in the enclosure of at least one of the eigenpair.
}
\label{Tab:performance}
\end{table}

\end{document}